\documentclass{amsart}
\usepackage{amssymb,latexsym}
\begin{document}

\title[ The Groups Aut(32$p$)]{The Automorphism Groups\\ of the Groups of 
Order $32p$}

\author{Elaine W. Becker}
\address{American Mathematical Society, Providence, Rhode Island}
\email{ewb@ams.org}
\author{Walter Becker}
\address{266 Brian Drive, Warwick, Rhode Island}
\email{w\_becker@hotmail.com}
\maketitle
\begin{abstract}

	The results of computer computations determining the
automorphism groups of the groups of order 32$p$ for $p \geq 3$ are
given in several tables. Presentations for the automorphism
groups of the groups of order 32, which in many cases appear
as direct product factors in the automorphism groups of order
$32p$, are also presented for completeness.

	Many of the groups of order  32$p$ with a normal sylow
$p$-subgroup have automorphism groups of the form:
Hol($C_p$)$ \times $Invariant Factor. A suggestion is made as to how one might
determine this invariant factor using only information on the automorphism group
of the 2-group associated with the group of order  32$p$, and the
normal subgroup of the 2-group associated with the extension
of the group of order  $32p$.

	Some general comments on the groups of order $32p^2$ and
their automorphism groups are made. A few explicit calculations
for the groups of order $32p^2$ are reported here. Knowing the automorphism 
groups for the groups of order $32p$ enables us to explicitly write down the 
automorphism groups for more than half of the automorphism groups of the 
groups of order $32p^2$.\\
\end{abstract}

\section{Introduction}

	In two previous papers we have discussed the groups and
the automorphism groups of the groups of orders  8$p$, 8$p^2$,  16$p$, and
16$p^2$ \cite{1,2}. In this paper we do a similar thing for the groups
of order 32$p$. The groups of order 32$p$ were determined many years ago by
Lunn and Senior \cite{3}. As part of his thesis on the construction of
solvable groups, Dr. R. Laue obtained presentations for the groups
of order 96 \cite{4}. The work of Miller (also cited in \cite{3}) and of 
Lunn and Senior in the 1930s do not give presentations for these groups,
and the relations given by Laue are rather more involved than they
need to be.\\

	As far as the authors know the automorphism groups for almost
all of these groups have yet to be determined. The object of this
report is to eliminate this gap in the known properties of these
groups. In many cases the automorphism groups for these groups
involve, as a direct factor, the automorphism groups of certain groups
of order 32. In view of the fact that the automorphism groups of the
groups of order 32 are also relatively unknown, we have included a
discussion of these groups here for the benefit of the reader. In
previous papers we have also included a discussion of the groups with
a normal subgroup of order $p^2$. We have not attempted a complete
description of the groups of order $32p^2$ here. Some general comments on the
groups and the automorphism groups of the groups of order $32p^2$ are
made in section 4.\\ 

      We began this work when we were graduate students at Michigan State 
University around 1980. One of the reasons for starting this work on 
automorphism groups of finite groups was the various conflicting opinions on 
just what was known about the automorphism groups of the groups of order 32. 
Some of these comments, we found puzzling were:
\begin{quotation}
The automorphism groups of the groups of order 32 are well known.
\end{quotation}
\begin{quotation}
The automorphism groups of the groups of order 32 can be found in the 
Hall-Senior Tables \cite{8}.
\end{quotation}
At the time we suspected that these automorphism groups were not all that 
well known. One reason that made us skeptical about just how well known the 
automorphism groups of the groups of order 32 were was derived from the 
fact that several of these automorphism groups have orders 32 or 64 and were 
not so identified as automorphism groups in the well-known reference work of 
Hall and Senior \cite{8}.

\section{The Groups of Order 32 and their Automorphism Groups}

The groups of order 32 were determined by Le Vavasseur \cite{5} and G. A. 
Miller \cite{6} around 1900. A minimal set of presentations for these groups 
was given by Sag and Wamsley in 1973 \cite{7}. Most of the relations used in 
this work are taken from the work of Sag and Wamsley. The eight 
presentations in Tables 1a to 1e, marked by [*] are different from those 
given in Sag and Wamsley's table. The groups of order 32 and their 
automorphism groups are given in Tables 1a to 1e. In many cases, a \textquotedblleft nice"  
description of just what these automorphism groups are is not known but 
presentations for these automorphism groups are given. In these tables the 
order structure is given for many of the groups and/or their automorphism 
groups. The following notation is used to specify the order structure of these
groups. For example, in the notes to Table 1a for the automorphism group of $C_8 
\times C_4$, we have\\ 

                            2-47-15  [$1^3,2^6,4^4,8^2$]. \\

This should be read as follows: there are 47 elements of order 2 distributed 
in 15 conjugacy classes. The elements are further broken down into three 
classes of one element each, six classes of two elements each, four classes of 
four elements each, and two classes of eight elements. The other cases 
follow in like manner. \\

A recent (October 2006) survey on the automorphism groups of finite 
$p$-groups is a paper by Geir T. Helleloid \cite{11}. This survey article contains 
many useful references dealing with various aspects of automorphism groups. 
We hope that our series of papers can at least serve as a useful adjunct to 
these articles. In particular we hope these papers can serve as a useful 
source of examples (and counterexamples) to illustrate various theoretical 
papers dealing with the automorphism groups of not only $p$-groups but more 
generally of finite, non-$p$-groups as well. \\

In regard to a comment by Helleloid \cite{11} on the automorphism groups of 
extra-special $p$-groups: \textquotedblleft Winter\ldots gives a nearly complete description of 
the automorphism groups of an extra-special $p$-group\ldots ". We have two of 
these groups here: numbers 42 and 43 of order 32. The automorphism groups of 
both of these groups are complete groups, a fact not mentioned by Winter or 
others as far as we know. Here we give in one case both the order and a 
presentation of this complete group of order 1920. In the other case, the 
automorphism group is just $S_4 \wr C_2$. This poses the question: Are the 
automorphism groups of other extra-special $p$-groups complete groups?

\section{The Groups of Order 32$p$ and their Automorphism Groups}

		The groups of order  32$p$ are given along with their automorphism groups 
in Tables 2a to 2e and Tables 3a to 3c. For those groups with a normal sylow 
$p$-subgroup the groups are listed in Tables 2a to 2e, and they are broken 
down according to the isoclinic class of the order 32 quotient group. The 
groups with a normal sylow 2-subgroup or groups without a normal sylow 
$p$-subgroup are listed in Tables 3a to 3c.

		The presentations of the groups given in Tables 2a to 2e were in many cases 
taken directly fom the computer input files, in order to reduce the 
possiblility of a transcription error.\footnote{This was true for the 
original manuscript written using a text editor. The presentations included 
below were re-edited, first using Microsoft Word and then converted at a 
later date to a LaTeX editor for this version.}

		The information given in Tables 2a to 2e not only allows one to construct 
the groups of order  32$p$ with a normal sylow $p$-subgroup but also gives 
their automorphism groups. A simple example should suffice to show how ones 
does this. Consider the cases in Table 2a for the 2-group ($2^2$,1), which is 
just
\begin{equation*}
				C_4 \times C_4 \times C_2.
\end{equation*}
	The presentation for this group is given in Table 1a:
\begin{equation*}
       a^4 = b^4 = c^2 = (a,b) = (a,c) =(b,c) = 1.
\end{equation*}

	This group can act on the cyclic group of order $p$ in two distinct ways, 
yielding two nonisomorphic groups of order  32$p$. Presentations for these 
two groups are:
\begin{gather*}
		a^4 = b^4 = c^2 = (a,b) = (a,c) = (b,c) = d^p =\\
		\begin{cases}(a^{-1})*d*a*d = (b,d) = (c,d) = 1 \qquad  \text{a 3}\\
             \textup{or} \\
			 (a,d) = (b,d)= (c^{-1})*d*c*d = 1. \qquad \text{c 4}
\end{cases}
\end{gather*}
	The entries in the column $C_2$ action in these tables specify the 
generator(s) of the order 32 groups, whose presentations are given in Tables 1a 
to 1e, that are acting on the cyclic group. The number next to the 2-group 
generator(s) gives the multiplicity of the group associated with the 
extension; see below for the details. The automorphism group is given by forming 
the direct product of the entry given in Tables 2a to 2e with the holomorph 
of the cyclic group of order $p$.

		The material on the automorphism groups of the groups of order  32$p$ is 
believed to be new. In those cases where the  32$p$ group has a normal sylow 
$p$-subgroup and is NOT of the form $C_p \times $ (a 2-group) (i.e., is not 
a direct product of its sylow $p$-subgroups), the automorphism groups usually 
(but not always!) take the form:\footnote{The exceptions here arise due to the
presence of an odd prime 
occurring in the invariant factor, e.g., cases \#1 and \#2 in Table 2a, Table 
2b cases 8a and 8c, as well as others in Tables 2a, 2b and 2d.}
\begin{equation*}
			\text{Hol}(C_p) \times \text {some 2-group}.
\end{equation*}
In some of these cases the 2-group has order 128 or greater and a unique 
description of the 2-group would be useful. In many cases they turn out to 
be direct products, e.g.,
\begin{equation*}
			[(C_2 \times C_2 ) \wr C_2 ] \times C_2 \times C_2 ,
\end{equation*}
but in other cases no such simple representation is available.\\

	Let $G$ = $C_p$  @ $X$,  where $X$ is a specific 2-group of order 32.
Let the action of $X$ on $C_p$ be by means of an operator of order 2.
Then we can view the action of $X$ on $C_p$ in $G$ as a mapping as
follows:
\begin{equation*}
			f: X  \mapsto \text{Aut}(C_p).
\end{equation*}
	The kernel of this mapping is called the \underline{group associated with 
the extension}. Let $H$ be the kernel of the mapping $f$. We conjecture that 
the invariant factors listed in Table 2 are those subgroups of the 
automorphism group of $X$ that are defined as holding the kernel of the 
homomorphism $f$ fixed. For each case in which the action of the 2-group on 
$C_p$ is of order 2, the group $H$ is given in the last column of Table 2 under
the heading \textquotedblleft group associated with extension".\\

		Table 4 shows the results of an attempt to verify this conjecture. The 
calculations listed in Table 4 were done as follows. The presentation for 
the group of order  32$p$ determines a fixed normal subgroup of the group of 
order 32, namely the group $H$. We asked CAYLEY to find all of those 
generators that CAYLEY obtained from getting the automorphism group of the group of 
order 32, and then to determine those generators from this set that mapped 
the group $H$ onto itself. We did not require that this mapping of $H$ onto 
itself be the identity mapping. In most cases this subgroup of the 
automorphism group of the order 32 group is isomorphic to the factor listed 
in Table 2. Denote this group by $T$. CAYLEY was then asked to see if the 
group $T$ was isomorphic to a normal subgroup of Aut($X$). In most cases this 
was in fact the case. The program used to do these calculations is described 
in Appendices I and II.\\

To simplify the discussion and show what was found let us consider a couple 
of specific examples.\\

\underline{Example 1.} Group number 48 of order 32. In this case there 
are two possible extensions: one associated with the group $Q_2 \times C_2$ 
(a characteristic subgroup of the order 32 group) and a second one $<2,2|2>$ 
with multiplicity 2 in that normal subgroup lattice. In the first case,  
Aut($G$) is just Hol$(C_p) \times $ Aut(32 \#48). In the second case the 
2-group factor is $D_4 \times D_4$. The group $D_4 \times D_4$ is a 
characteristic subgroup of Aut(32 \#48).\\

Example 1 illustrates a rather common pattern found in those cases in which 
there is only one extension with a 2-group factor of half the order of the 
order 32 groups' automorphism group.\\

\underline{Example 2.} In other cases, e.g., group 16 of order 32, we have 
more than one possible extension which gives rise to a 2-group of order 128. 
Two different extensions have isomorphic automorphism groups and the 
2-factor is not a characteristic subgroup of Aut(32 \#16) but rather 
appears with multiplicity 2 in that normal subgroup lattice.\\

In other cases, e.g, for $[C_4 \text{Y} Q_2] \times C_2$ the 256 factor is 
the sylow subgroup of Aut$([C_4 \text{Y} Q_2] \times C_2$) and is not a 
normal subgroup. In these cases it is harder to verify that one does get a 
subgroup of the automorphism group of the 2-group. This will get worse as 
the order of the 2-group increases. This is particularly true in those cases 
for which the calculations described below give the wrong order group.\\

		In a few cases noted in Table 4 the order of $T$ was not equal to the 
order of the invariant factor in Table 2. In those cases that we have 
checked (numbers 8c, 16ab, 19a) as well as for numbers 6, 15, 17 and 41 
which were not run, the given factors all are subgroups (but not necessarily 
normal subgroups) of the automorphism group of the order 32 quotient group 
of this group of order  32$p$. For cases 2b and 8c, no attempts were made to 
see if the associated group given in Table 2a is a subgroup of the 
automorphism group of the order 32 group.\\

In view of the number of correct identifications of $T$ with the factors in 
Table 2, this conjecture seems to be a close approximation to what is actually 
happening here. The problem here may be that the method used to calculate 
the group $T$ in Table 4 is incorrect, or that the conjecture needs to be 
modified. Problem: Why do some of the calculations give the wrong group? The 
wrong group arising in this way has, so far, always been of the wrong order 
so that the incorrectness of the calculation is easy to verify. This may or 
may not always be the case, so why this is happening would be useful
to know. Even more useful would be to have a method to calculate these 
automorphism group invariants directly from the automorphism group of  the 
2-group.

\section{Some General Comments on the Groups and Automorphism Groups of
the Groups of Order $32p^2$}

When we started to construct the groups of order $32p$ and to compute their 
automorphism groups (circa 1980) the number of groups of order $32p^2$ was 
unknown, at least to us. As can be inferred from the situation occurring in 
the case for the groups of order 16$p$ and $16p^2$ the number of groups goes 
up very rapidly as one increases the order of the normal $p$-subgroups, for 
a given 2-group. When this work was begun, the work of H. Besche and B. Eick 
\cite{12} was in the future. This work of Eick and Besche has determined 
the number of groups of order $32p^2$ for the primes $p=3$ (1045), 5 (1211) 
and 7 (928).\\

These are the only orders in this sequence of groups of order $32p^2$ that 
are covered by the work of Besche and Eick. It might be instructive to 
include here how we obtained some information on the possible number of 
groups of these orders, without having recourse to the work of Eick and 
Besche. At that time we had a number of ways to get a rough estimate of the 
number of groups one might find in various orders. One crude estimate is to 
look at the number of groups of order $2^n$ and of order $2^{n-1}*p$. These 
are roughly the same order of magnitude. A crude estimate of the number of 
groups one might expect to find is obtained by observing the number of groups 
found in the orders $2^n$, $2^{n-1}*p$ and $2^{n-2}*p^2$ for some values of 
$n$. For simplicity we look only at the case for $p=3$.\\

\begin{tabular}{|c|c|c|c|c|c|} \hline
\multicolumn{6}{|c|}{The number of groups of order $2^k*p^m$ }\\ \hline
$|G|$ & $n=4$& $n=5$&$n=6$&$n=7$&$n=8$ \\ \hline
$2^n$  &  14 & 51 & 267 & 2328 & 56,092 \\ \hline
$2^{n-1}*p$ & 15 & 52 & 231 & 1543& 20,169 \\ \hline
$2^{n-2}*p^2$& 14 & 50 & 197 & 1045 &8681  \\ \hline
\end{tabular}\\

One observes here that very roughly the number of groups of order $2^n$ and 
the number of groups of order $2^{n-1}*3$ are of the same order of 
magnitude. That this is not a very good or accurate estimate can be seen by 
looking at the now known results for orders 64, $32p$ and $16p^2$ (for the 
case of $p=3$). One could do a little better by looking at cases where $p=5$ 
or 17, but this did give us a crude estimate of what would be involved if we 
went on to try to determine the groups of order $32p^2$ and their 
automorphism groups. Note also the number of groups of order 192 and 256 
were also not known when we started this work in the early
1980s.\footnote{The calculations reported in Tables 1, 2 and 3 were completed by 
the late 1980s and an early version of this paper was written in the 
early 1990s.} \\

Another way to get an estimate of the number of these groups is to look at 
the Hall-Senior Tables, or more accurately the charts of the lattice of 
normal subgroups contained therein. This is how we got a more accurate 
estimate on the number of groups of order 192. One first must realize that 
the number of groups of order $2^n*p$ where the 2-group acts on the 
$p$-group is determined by how many inequivalent normal subgroups the 
2-group has. This is shown in the Hall-Senior Lattice charts. The number of 
semi-direct products $C_p @ G[2^n]$ where the 2-group action is by an 
operator of order 2 (i.e., $C_2$) can be found by just adding up the number 
of normal subgroups of order $2^{n-1}$ for the groups of order $2^n$. For 
the order 32 groups one finds there are 144 such cases. In the older 
literature, e.g., G. A. Miller and others, also Lunn and Senior (\cite{3}, pp.  
322 and 324), these are called dimidiations. For example in \cite{3} one 
finds the statement \textquotedblleft for every prime $p$ the number of groups of this 
division is equal to the number of distinct dimidiations of the groups of 
order 32. This number has been shown to be 144." The orders of the 
automorphism groups can also be found from these charts. For the case when the 
normal subgroup associated with the extension is a characteristic 
subgroup, the automorphism group is just Hol$(C_p) \times$ Aut(G[32]), where 
G[32] is the group of order 32 in question. If the group associated with the 
extension is not a characteristic subgroup, then the order of this $32p$ 
group's automorphism group is just $p*(p-1)*|\text{Aut(G[32])}|/k$, where $k$ 
is the \textquotedblleft multiplicity" of this subgroup which is associated with the 
extension. These numbers are also given in these subgroup lattice charts. 
This is how we guessed that the number of groups of order 192 would be quite 
large, i.e., 267 direct products plus at least 1120 cases in which the 
action by the 2-group on the $p$-group is by means of an operator of order 
2. If the reader wants to check this, one can just take the Hall-Senior 
charts and start counting them. This, by the way, will tell you just how many 
groups of order 192 you will get in this way from a given group of order 64. 
It might also be of interest to note that 700 of these 1120 groups have a 
characteristic subgroup for the group associated with the extension, and 
hence from above one knows what these automorphism groups are, i.e., 
Hol$(C_p) \times$ Aut([G[64]) (see \cite{13}).\\

For the number of $C_4$ action cases, knowing the number of such cases in 
order $32p$ (which is 40 when $p=5$), we know that we will have 200 cases in 
order $32*25$. These arise as follows. There are 40 from the case when the $p$-group 
is $C_{25}$ corresponding directly to those in order $32*5$ arising from a 
$C_4$ action on $C_5$. We get another $4\times 40$ from the $C_5\times C_5$ 
cases for a total of 200. The corresponding actions here (for the 
$C_5\times C_5$ cases) are:
\begin{equation*}
a^5=b^4=a^ba^2 =c^5 = (a,c)= \begin{cases} (b,c) \\ c^bc \\ c^bc^2 \\ 
c^bc^3\\ \end{cases}=1.
\end{equation*}
Each of these relations will give rise to a distinct (nonisomorphic) group 
of order 800. Note here for the  $p=3$ case, i.e., $(C_3\times 
C_3)$@(2-group) with a $C_4$ action we have 40 cases here and these 
automorphism groups are $[144] \times $ invariant factor. This [144] factor 
here is the complete group $(C_3\times C_3) @ QD_8$. For the $p=5$ cases 
this [144] factor goes over into Hol$(C_{p^2})$ for the case when the 
$p$-group is $C_{p^2}$, and for the other cases as follows:
\begin{equation*}
a^5=b^4=a^ba^2 =c^5 = (a,c)= \begin{cases} (b,c)=1\rightarrow 
\text{Hol}(C_5) \times C_{4}\\
c^bc=1\rightarrow \text{Hol}(C_5) \times \text{Hol}(C_5) \\ c^bc^2 =1
\rightarrow \text{Hol}(C_5 \times C_5) \\ c^bc^3 =1\rightarrow 
\text{Hol}(C_5 ) \wr C_2.\\ \end{cases}
\end{equation*}
The generalization to other primes with $p\equiv 1$ mod(4) is
straightforward; just replace 5 by $p$ and 4 by $p-1$ in these automorphism 
groups.\\

One can also continue as in this fashion for the $C_q$ group actions for 
larger 2-groups, e.g., $C_8$, $C_{16}$, etc. As one goes up in order for the 
2-group one gets more and more \textquotedblleft quotient groups" becoming
involved and the 
number count can be quite daunting using charts such as that found in \cite{8}. 
At some point this method really does not become very useful for
\textquotedblleft global" 
aspects of groups of order $2^np^2$ (or for that matter for groups of order 
$2^np^m$). This method may still prove to be a very useful guide if one is only 
concerned with a single 2-group, or very few 2-groups. This will tell one 
just how many extensions of a particular type to expect, and yield some 
information on what their automorphism groups might be.\\

The cases for which the action is by means of a group of order 16 or 32 is 
somewhat more involved, since these cases arise only in special orders, 
e.g., $QD_8$ acts on $C_p \times C_p$ only when $p \equiv  1$ or 3 mod(8), 
and $C_4 \wr C_2$ acts only when $p \equiv  1$ mod(4). The material in 
Table 5 shows what may be expected for various order 16 actions arising from 
different groups of order 32. (For the final word on this matter one should 
consult the work of Besche and Eick mentioned above.)\\

For the case of groups of order 32, see also the appendix to the authors' 
paper on groups of order $16p^2$, where subgroups of the sylow 2-subgroups 
of Aut($C_p \times C_p)\simeq$ GL(2,$p$) are given for a variety of 
primes.\\

For the other cases, such as that with a $C_2\times C_2$ action, however, a 
different method appears to work for us. Each group of order $8p^2$, $16p^2$ 
and $32p^2$ that is not a direct product of its sylow subgroups has a 
nontrivial normal subgroup associated with it called the \textquotedblleft group asociated 
with the extension". In the case when this \textquotedblleft group associated with the 
extension" has a quotient group that is a direct product, e.g., $C_2\times 
C_2$ or $C_4\times C_2$, \ldots , we can have two or more groups of order $8p^2$ 
or $16p^2$, $32p^2$ \ldots\ possessing this specific group as the \textquotedblleft group 
associated with the extension". This appears not to be the case for the 
normal subgroups whose quotient groups are $C_2$, $C_4$, $C_8$, $D_4$, 
$Q_2$,\ldots , i.e., when the quotient group is NOT a direct product. The number 
of groups arising in order $2^np^2$ due to a $C_2$, $C_4$, $C_8$, \ldots\ action 
on the $p$-group is just equal to the number of such quotient groups in the 
2-group's lattice. This observation is probably the result of some well-known
group theory argument, but it is unknown to us at this time.\\

The idea here is first to just count the number of lines emanating from the 
2-group's normal subgroup of order $2^{n-1}$ and ending up on a second 
normal subgroup of order $2^{n-2}$, whose quotient group is $C_2\times C_2$. 
This information is also contained in the Hall-Senior charts. If one does 
this for the groups of order 8 and order 16, one gets the numbers 6 and 35, 
respectively. The number of groups of order $8p^2$ arising from a $C_2\times 
C_2$ action  is 6, and the corresponding number for the groups of order 
$16p^2$ is 35. If this correspondence holds for higher-order 2-groups, e.g., 
groups of order 32 or 64, then the number, as well as the distribution of 
such cases among the 2-groups of the appropriate order, can be \textquotedblleft easily" (?) found from the lattice diagrams in the Hall-Senior Tables \cite{8}. Using 
this argument one gets 263 cases. The number of such groups found by Eick 
and Besche is 274. \\

Continuing on in this vein one might get an estimate on the number of groups 
with a $C_4\times C_2$ action to be found in much the same way. That is, for 
those normal subgroups whose quotient group is $C_4\times C_2$, or indeed 
some other direct product, one might expect that, for this given normal 
subgroup (of order 8 here) that the number of groups of order $32p^2$ with 
this normal subgroup being the \textquotedblleft group associated with the extension"
is equal to the number of cases of lines incident on this normal subgroup 
originating from a normal subgroup of order 2 or larger. A count of the cases 
for $C_4\times C_2$ gives us 91 cases incident upon 41 $C_4\times C_2$ 
quotient groups. The Eick-Besche table gives us 112 cases for $p=5$. \\

These considerations for the $C_2\times C_2$ and $C_4\times C_2$ cases do give us 
a rough idea of the number of such extensions to be expected. It would 
thus appear that either our conjecture about the number of such groups being 
equal to the number of incident lines is wrong or that some lines are 
missing in the original Hall-Senior charts. It would be nice to know which 
is correct. If the problem is with missing lines in the Hall-Senior charts, 
it would be useful to have an errata for the appropriate group lattice 
diagrams listed somewhere for interested readers. This may already be the 
case in some of the group packages distributed with GAP (and MAGMA ?). 
However, the lattice diagrams in the GAP case require a Unix or Linux 
operating system, which might limit access to this material for those who use other 
operating systems on their computers.\\

Knowing the automorphism groups for the groups of order $32p$ when the 
action of the 2-group on $C_p$ is by $C_2$, we can also state what the 
automorphism groups are for those groups of order $32p^2$ when the 
action of the 2-group on the $p$-groups is by the group $C_2$. The number of 
groups of order $32p^2$ where the $p$-group is $C_{p^2}$ in this case is 
144. These groups' automorphism groups can be written down by looking at the 
corresponding case for order $32p$ and replacing the group 
Hol$(C_p)$ by Hol$(C_{p^2})$. For the case when the $p$-group is $C_p\times 
C_p$, the 2-group can act on the $p$-group in two different ways. It can act 
as an operator of order 2 on the first $C_p$ and commute with the second 
$C_p$.  This gives us another set of 144 groups of order $32p^2$. These 
groups' automorphism groups are obtained from the corresponding order $32p$ 
case by replacing the Hol$(C_p)$ by Hol$(C_p) \times C_{p-1}$. In the second 
case, the 2-group acts as an operator of order 2 on both $C_p$ factors, 
yielding a third set of 144 groups of order $32p^2$. These groups' 
automorphism groups are obtained by replacing the Hol$(C_p)$ group in the 
corresponding $32p$ case by the group Hol$(C_p\times C_p)$. Hence at this 
point we know the groups and automorphism groups of the $2\times 51 $ direct 
products as well as $3 \times 144$ groups = 504 groups of order $32p^2$. \\

From our previous results \cite{1,2}, we expect that many of the 
automorphism groups arising from groups of order $32p^2$ with a $C_2\times 
C_2$ action will take the following form:
\begin{quotation}
1. $\text{Hol}(C_p)\times \text{Hol}(C_p) \times$ an invariant factor from 
Table 2.
\end{quotation}
Note here that in a few cases, however, a new lower-order invariant factor not 
present in the $C_2$ invariants listed in Table 2 will appear:
\begin{quotation}
2. $ \left[ \text{Hol}(C_p)\times \text{Hol}(C_p) \times \text{ invariant 
factor}\right]@C_2$.
\end{quotation}
The $C_2$ in this second form acts on $\text{Hol}(C_p)\times 
\text{Hol}(C_p)$ so as to yield the wreath product Hol$(C_p) \wr C_2$. In 
most cases when the 2-group has more than one $C_2$ action on $C_p$ and 
yields invariants whose orders differ by a factor of two, the action of the 
$C_2$ on the smaller \textquotedblleft invariant factor" yields the larger
$C_2$-invariant factor. The form (1) occurs in the \textquotedblleft direct action" cases, i.e., when each 
one of these two $C_2$'s acts on one of the $C_p$'s and commutes with the other 
$C_p$ and the second $C_2$ acts on the second $C_p$ and commutes with the 
first $C_p$. The second case arises in the \textquotedblleft cross action" cases, i.e., when 
one of the $C_2$'s acts on both $C_p$'s and the second $C_2$ acts on only 
one of the two $C_p$'s. Many examples of these types will be displayed in a 
forthcoming report \cite{14}.\\

To determine the other cases, i.e., when we have a normal sylow 2-subgroup 
and cases without any normal sylow subgroup, one can proceed as follows. If we 
have the number of groups of order $2^np$ with a normal 2-sylow subgroup, we 
can get a good idea of just how many we have in the case of order $2^np^2$. Namely,
take all of the groups with a normal sylow 2-subgroup 
(\underline{omitting here those groups that are the} \underline{direct 
product of their sylow subgroups}) and multiply them by 2. This should give 
a good approximation to the number of groups with a normal sylow 2-subgroup 
in order $2^np^2$. The cases that are omitted here are those with an 
automorphism of the form $C_p\times C_p$, which are usually few in number. 
In the case of the groups of order 32, for example, we have just three cases, 
namely the elementary abelian group of order 32, $Q_2\times C_2\times C_2$, 
and $Q_2$Y$Q_2$ (number 42 in the Hall-Senior Tables). For the groups of order 
$32p^2$ (for $p=3$) these cases will give us 9 groups (three groups from 
each of these three groups) with a normal sylow 2-subgroup; all of the 
others, with a single factor of 3 in their automorphism group order will 
give just two cases, one when the $C_3$ action is from one of the $C_3$'s in 
the group $C_3\times C_3$ and one from a $C_9$ group acting as an operator 
of order 3 on the order 32 group. For the case of order $32p^2$ one finds 41 
such groups (for $p=3$). Note, however, if you ask the program GAP (and 
presumably CAYLEY or MAGMA) to give you the number of groups with a 
normal sylow 2 (or a normal sylow $p$) subgroup, it will return 41 + (2 times 
51) = 143 groups (or for the case of a normal sylow 3-subgroup, we get 914 
instead of 812). That is, it will also count the number of such cases 
arising from the direct product of their sylow subgroups in this tabulation, 
which may not be what the user wants to look at.  The automorphism groups 
may also be inferred as follows. For every case of a normal sylow 2-subgroup 
(again omitting the direct product of sylow subgroups cases) of the form 
(2-group)@$C_3$, the automorphism group of the the corresponding case of 
$2^np^2$, i.e., ([2-group]@$C_3$) $\times C_3$ and [2-group] @ $C_9$ acting as 
the $C_3$ in the (2-group) @ $C_3$ case, the automorphism groups are:
\begin{equation*}
\text{Aut([32]} @ C_3) \rightarrow \begin{cases} \text{Aut([32]} @ C_3) 
\times S_3 \quad (C_3\times C_3 \text{ case}), \\
            \text{Aut([32]}@C_3) \times C_3 \quad (C_9 \text{ case}). 
\end{cases}
\end{equation*}
The other three cases of groups with a normal sylow 2-subgroup arising in 
order $32p^2$ for $p=3$ come from the groups:
\begin{align*}
&A_4 \times A_4 \times C_2, \qquad \text{i.e., the extension} \quad 
\left(C_2\times C_2\right)@C_3 \times  \left(C_2\times C_2\right)@C_3 \times 
  C_2, \\
&Q_2@C_3 \times \left(C_2\times C_2\right)@C_3\simeq SL(2,3)\times A_4, \\
&Q_2\text{Y}Q_2@\left(C_3\times C_3\right).
\end{align*}
The automorphism groups of these groups are:
\begin{gather*}
A_4\times A_4\times C_2\rightarrow S_4\wr C_2, \\
Q_2 @C_3\times \left(C_2\times C_2\right)@C_3 \rightarrow S_4\times 
S_4\times S_3, \\
Q_2\text{Y}Q_2@\left(C_3\times C_3\right)\rightarrow S_4\wr C_2.
\end{gather*}
A systematic way of looking for the groups without a normal sylow subgroup will be 
discussed in connection with the groups of order $64p$ without a normal 
sylow $p$-subgroup. Given the set of groups of order $2^np$ without a normal 
sylow subgroup one can proceed as follows. One can use an argument similar 
to that for the normal sylow 2-sugroup case here. Every \textquotedblleft nonnormal
sylow group" of order $2^np$ is also a \textquotedblleft nonnormal sylow group"
case for order 
$2^np^2$ when we form the direct product of the $2^np$ group with the group 
$C_p$. Likewise one can expect to get a second set of nonnormal sylow types 
by replacing the $C_p$ group in the $2^np$ cases with the group $C_{p^2}$ 
with the group $C_{p^2}$ acting as an order $p$ operator (i.e., as a $C_p$) 
on the 2-group in question. Additional cases can arise from groups of order 
$2^np$ with a normal sylow 2-subgroup, but which possess a direct factor such as 
a $C_2$, or some other 2-group which can act directly on the second $C_p$ 
here. Consider the case of the group $\left(C_2\times C_2\times C_2\times 
C_2\right)@C_3 \times C_2$. Now add a second $C_3$ and we have 
$\left(C_2\times C_2\times C_2\times C_2\right)@C_3 \times C_3@C_2$, which 
is a group without a normal sylow subgroup.

This turns the previously normal sylow 2-subgroup case into a 
nonnormal sylow subgroup case. This is a case-by-case (or group-by-group) 
approach and again can get rather tedius rapidly as one increases the orders 
of the groups involved. The main advantage here is if one is dealing with a 
specific group or a small set of groups of interest in which case this approach 
may yield addional information of interest to the user. One should also 
point out here that these cases in which we have groups of order 
$2^np^m$ with normal sylow 2-subgroups or no subgroups with a sylow subgroup 
arise only for a discrete set of primes $p$, and usually small ones at that. 
As the prime $p$ increases, the number of such cases drops rapidly and 
eventually is zero for sufficiently large primes. \\
\begin{table}
\begin{center}
\begin{tabular}{|c|c|c|c|c|c|c|c|c|c|c|}\hline
\multicolumn{11}{|c|}{Summary for the groups of order $32p^2\,\,\dagger$ } 
\\ \hline
    & \multicolumn{8}{|c|}{cases with a normal sylow subgroup} & & \\\hline
    & \multicolumn{7}{|c|}{image of 2-group $\ddagger $} & \# sylow 2- & no 
sylow&direct \\ \cline{2-8}
$p$ & $ 1$&$ 1^2$&$ 2$& $3$&$D_4$&$Q_2$& $16$ & subgroups&subgroups&products\\ 
\hline
3 & 432&274&40&9&42&11&4&41&90&102\\
5*&432&274&200&9&42&11&20&4&4&102\\
7*&432&274&40&9&42&11&11&2&4&102\\ \hline \hline
Aut& yes& &yes&& see& see & & yes&some but & yes\\
known   & &&& &\cite{9}&\cite{9}&&&not all&\\   \hline
\multicolumn{11}{|l|}{$\ddagger $ $1\equiv C_2$ image, $1^2\equiv C_2\times 
C_2$ image, $2\equiv C_4$ image, $3\equiv C_8$ image,}\\
\multicolumn{11}{|c|}{and 16 $\equiv$ groups of order 16}\\ \hline
\multicolumn{11}{|l|}{* For $p=5$ we have 112 cases from a $C_4\times C_2$ 
action and} \\
\multicolumn{11}{|c|}{one from an order 32 image.}\\
\multicolumn{11}{|l|}{* For $p=7$ we have one case with an order 32 image.}\\ 
\hline
\multicolumn{11}{|l|}{$\dagger$ Based upon results in Small Group Library of 
Eick et al. \cite{12}} \\ \hline
\end{tabular}
\end{center}
\end{table}

\section{  Conclusions and Future Considerations}

	In our first three papers we have given explicitly the automorphism groups 
for the groups of orders  8$p$, $8p^2$,  16$p$, $16p^2$, and  32$p$. For the 
most part this work has been in the nature of just reporting upon the 
results derived from a large scale computer study of these finite groups and 
their automorphism groups. Certain systematic behavioral patterns for these 
automorphism groups have been pointed out, notably the relatively simple 
dependence of the automorphism groups upon the $p$-groups. In the cases 
studied so far we have held the 2-group fixed and determined how the 
automorphism groups changed as we changed the $p$-group. The next major step 
in a systematic study of the automorphism groups would be to see how the 
automorphism groups change when one varies the 2-group. There are a number 
of ways to interpret the question:\\

		 See how the automorphism groups of finite groups of orders $2^n * p$ 
and/or $2^n * p^2$ vary as a function of the 2-group (i.e., as a function of 
$n$).\\

	One such way would be to look at extensions of groups of the form $G$[$p$] @ 
$A$, where $G$[$p$] is either $C_p$ or $(C_p \times C_p)$ and $A$ is a 2-group, and ask 
how are the automorphism groups of these groups related to the automorphism 
groups of the groups $G$[$p$] @ $B$, where the 2-group $A$ is a normal subgroup of 
index 2 in the 2-group $B$. In many cases this is equivalent to asking how are 
the automorphism groups of the 2-groups $A$ and $B$ related, which is again a 
question that we have not addressed in these papers. The large numbers of 
explicit examples worked out in these papers may provide clues as to what are 
the appropriate theorems to look for.\\

	Another way to look at this same question is to see how the automorphism 
groups change with the 2-group, but also requiring that the action be by the 
same group, e.g., a $D_4$ or a $Q_2$. This is the approach that will be 
reported on in \cite{9}.\\

	Some other questions that arise in this study were already
mentioned above:\\
\begin{quotation}

a.) Can one in fact construct those automorphism group invariants listed in 
Table 2 directly from the knowledge of the group of order $2^n*p$ and the 
automorphism group of the related 2-group?
\begin{quotation}
The results given in Table 4 imply this might be the case. The discussion 
above however did point out some problems with a too simplistic way of 
approaching the calculation of these factors. If this idea can be pushed 
through, it would greatly simplify the calculation of large classes of 
automorphism groups of finite order. This method would probably be mainly 
useful for actions involving the cyclic groups. The $p$-dependence for the 
other groups seems a bit more complicated to handle in this manner.
\end{quotation}
b.) It would be nice to have a simple explanation for the observation 
relating the patterns found in the Hall-Senior Tables' lattice diagrams with 
the number of finite groups of the form $p^m *2^n$, for $n$ less than or equal 
to 5 or 6.
\begin{quotation}
		In the analysis of the groups of order $16p^2$ in which the action of the 
2-groups on the $p$-group was $C_2 \times C_2$, one of the authors had the 
feeling that one could just look at the $(C_2 \times C_2)$ actions on the 
$p$-group and tell if their automorphism groups would be isomorphic. To be 
more precise, the conjecture is: all of the groups of order $16p^2$ (arising 
from a $C_2 \times C_2$ action) with the same automorphism group could be 
made to have the same $C_2 \times C_2$ actions on the $p$-group. One should 
note here that in these cases the 2-groups are not isomorphic. The number of 
cases dealt with there was comparatively few in number to what is supposed 
to appear in the orders $32p^2$. Therefore the order $32p^2$ groups might be 
a better place to test out this conjecture if one cannot come up with a 
more traditional group-theoretic
proof of the correctness of this assertion or show that this 
conjecture is false.
\end{quotation}
\end{quotation}

\section{Acknowledgements}

	This work was done over a period of years at several different 
institutions. We started this work while we were at Michigan State 
University. The bulk of the $p=3$ and $p=5$ cases, however, was done at the 
University of Rhode Island and at Syracuse University in the mid to late 
1980s. The work was finished up later on at Brown University. The 
calculations using the minimal presentations of Sag and Wamsley were all 
redone at Brown University, using the computers in the Department of 
Cognitive and Linguistic Sciences, in order to simplify the method of 
presentation. The acknowledgements mentioned in the previous paper on groups 
of orders  8$p$ and $8p^2$ are also relevant not only here but throughout 
this entire series of papers. The major difficulty the authors have had over 
this long span of time is the lack of personal contact with group theorists and 
others involved with computational group theory which would have made the 
effort more interesting (and enjoyable) as well as improving the 
presentations given in these papers. The authors would appreciate any 
comments or suggestions on the mode of presentation or requests for 
inclusion of additional information on these groups
or their automorphism groups.\\

\newpage

\section{Appendix I. Outline of Autsubc:\\ The CAYLEY Program for Table 4. }

\subsection{General orientation}
In this appendix we give a discussion of the program that was used to 
calculate the entries given in Table 4. Initially we present a fragmented 
version of the program, detailing what each part did. Later on we just list 
the entire program. One should point out here that we used CAYLEY, which used 
a \textquotedblleft determinative" rather than a random method for calculating various 
properties of the input group.  The program as presented below needed  an 
input from a separate program. The program below needed to have a specified 
normal subgroup of $G$ [in this case for Table 4, a specific subgroup of a 
group of order 32]. This was determined in a separate run from another 
program. We shall outline this program and then proceed on with a discussion 
of Autsubc in Appendix II. This way of doing things required us to make sure 
that every time we specified the same group $G$ as our input, the ordering of 
the output (here the normal subgroups of $G$) would appear in the same order or
sequence. In GAP, at least for some properties of a group, calculated with 
GAP and returned as a list, e.g., the generators of the automorphism group 
of $G$, this is not the case. In fact, the list of generators may not even be 
the same list, permuted about, but may contain different elements from run 
to run. In view of this, most GAP sessions should be run without 
interruption, rather than basing the input of a second run on the output of a 
previous run.\\

The routine Autsubc needs a group called nx for its input. This group is the 
\underline{group associated with the extension} for a specific group of 
order $96 \simeq C_3@[32]$. The way this group was determined is given in 
the following schematic CAYLEY program.\\
\begin{quotation}
We start out with
\end{quotation}
g:free(a,b,c);\\
g.relations:$a^{16}=b^2=(a,b)=c^3=a^c*a=(b,c)=1$;
\begin{quotation}
This specifies the group of order 96. The next set of instructions gets the 
normal subgroup of interest, nx.\\
\end{quotation}
p=3;\\
n= normal subgroups(g);\\
L=length(n);\\
sp=sylow subgroup(g,p); \\
cp=centralizer(g,cp);
\begin{quotation}
This determines the normal subgroup of $G$ that commutes with $C_3$.
\end{quotation}
scp=sylow subgroup(cp,p);
\begin{quotation}
This is just the group $C_3$.\\
\end{quotation}
qcp=cp/scp;
\begin{quotation}
This is the factor group cp by scp, which is the group associated with the 
extension of $G$.
\end{quotation}
for i=2, L do\\
    x=n[i];\\
    if(qcp eq x) then\\
      print i,x;\\
     end;\\
end;\\
nx=n[i];\\

The group nx is what is used to proceed with Autsubc.\\

The first part of Autsubc does the following. For a given group $G$ [here of 
order 32] this fragment takes a specific normal subgroup of $G$ [here a group 
of order 16], called nx, and finds the generators of Aut($G$) that map the 
subgroup nx back onto itself. It then finds the subgroup, h, of Aut($G$) that 
this subset of generators of Aut($G$) generates. The conjecture, or hoped for 
result, is that h is isomorphic to the invariant factor in Aut[$C_3 @ G] 
\simeq S_3\times$ inv.  This takes us down to part f, where we have XZ = x 
eq h appearing in the program. Here x is a normal subgroup of Aut($G$) that 
can be expressed as n[some number $m$] that can be compared with the initial 
input nx = n[some number $m'$] to see if the two arguments of n are 
equal.\\

\subsection{Outline of program.}
The following is a general outline of the program Autsubc:\\
a. Given $G$[32] we find
\begin{quotation}
   1. the normal subgroups of $G$,\\
   2. the automorphism group of $G$.\\
\end{quotation}
b. We input here nx = n[number from above run], the invariant factor, and find 
its conjugacy classes and automorphism group to check on the correctness of 
this input group.\\
c. We find the generators of Aut($G$) that map the normal subgroup nx in $G$ 
back onto itself. The subgroup of Aut($G$) that these generators generate is 
called h and we determine its order.\\
d. If the order of nx and h are the same, then we set the command
\begin{equation}
\text{ngiven2 =ngiven};
\end{equation}
otherwise, we set
\begin{equation}
\text{ngiven2=nng2}.
\end{equation}
To properly make use of the (nng eq 1) command, one needs to run the program 
twice, once to get the group h and a second time to continue properly from 
here. One could rewrite the program to do this in one pass, which if one 
uses GAP one might need to do.  \\
e. Compute properties of the group h, determined as a subgroup of Aut($G$):
\begin{quotation}
1. Is h a normal subgroup of Aut($G$)?\\
2. Is h a direct product of two groups?\\
3. Get the conjugacy classes, order structure and inner automorphism group 
of h.
\end{quotation}
\newpage
\noindent f. Properties of Aut($G$):
\begin{quotation}
1. Calulate Aut( Aut($G$) ). \\
2. Find all normal subgroups of Aut($G$) whose order and number of conjugacy 
classes equals those of the group h, and see if one of these normal 
subgroups is equal to the group h.\\
3. For all normal subgroups with $|$x$|=|$h$|$ and ncl(x)=ncl(h), see if x is a 
direct product.
\end{quotation}
g. Mapping of subgroups of Aut($G$) by Aut( Aut($G$) ):
\begin{quotation}
1. Determine the characteristic subgroups of Aut($G$).\\
2. Determine how the noncharacteristic normal subgroups of Aut($G$) are 
mapped about by Aut( Aut($G$) ).
\end{quotation}

\section{Appendix II. The Program Autsubc:\\ The CAYLEY Program for Table 4.}

The program for determining subgroups of Aut($G$) when a specified normal 
subgroup, nx[xy], of $G$ is held fixed.\\

The comments that appeared in the CAYLEY program are printed below in 
italic lettering.\\

\subsection {a. The initial input.}

The first fragment serves as a check on the specification of the group nx 
from the previous run, i.e., in the text below, the 12 in the line nx=n[12].\\
{\itshape \textquotedblleft This is for groups of order $p^6*q$."}\\
{\itshape \textquotedblleft Determine the characteristic subgroups of the group G and return them as the sequence chars."\\
\textquotedblleft This is number $C_{16}\times C_2$  b."}\\
p=2;\\
q=3;\\
{\itshape \textquotedblleft If ngiven2 eq ngiven, then mappings among normal subgroups only check out those for orders of the required normal subgroup."\\
\textquotedblleft If use ngiven2 eq 2, then it does the entire normal subgroup set of 
mappings."\\
\textquotedblleft If use ngiven2 different than ngiven, then set nng=2."\\
\textquotedblleft The starting point for mappings is then for order nng2."}\\
nng2=2;\\
nng=1;\\
n4=4; {\itshape \textquotedblleft Search for direct product with first factor of 
order n4 or more."}\\
chars=empty;\\
G:free(a,b);\\
G.relations:$a^{16}=b^2=(a,b)=1$;\\
pring g, order(g);\\
n=normal subgroups(g);\\
L=length(n);\\
print L;\\
autg=automorphism group (g);\\
print order(autg);\\
chars=empty;\\
nautgens=ngenerators(autg);\\
{\itshape \textquotedblleft Determine which normal subgroups of G are left invariant by the automorphism group."}\\
nx=n[12]; {\itshape \textquotedblleft Specify here the group associated with the 
extension."}\\
orf=16; {\itshape \textquotedblleft orf is order of aut(g) factor being looked for."}\\

\subsection{b. Check on the input group nx.} This part just consists of the 
command\\
print nx,classes(nx). \\

\subsection{c. Determination of the group h, a subgroup of Aut($G$)} The 
first part of the routine here determines the subgroup of Aut($G$) that maps 
the normal subgroup nx back onto itself. This is done by determining those 
generators of Aut($G$) that map nx onto nx and then taking these automorphisms 
of Aut($G$) and getting the subgroup of Aut($G$) that these generators 
generate.\\
ux=automorphism group(nx);\\
print order(ux);\\
L=ngenerators(autg);\\
for i=1 to L do
\begin{itemize}
  \item{f=automorphism(g,autg.i);}
  \item{invar=f(nx) eq nx;}
  \item{if(invar) then}
     \item{chars=append(chars,autg.i);}
   \item{end;}
\end{itemize}
end;\\
print length(chars);\\
h=$<$chars$>$;\\
print order(h);\\
ngiven=order((h);\\
\subsection{d. Check on the orders of nx and h.} The following few lines 
determine how we handle the normal subgroups of Aut($G$). If nng equals one, 
then we proceed to calculate various properties of the normal subgroup h, 
such as its class/order structure, is it a direct product, etc. If nnq 
differs from one, then we just list the normal subgroups of Aut($G$) from the 
next part of the program.\\
if(nng eq 1) then
\begin{itemize}
   \item{ngiven2=ngiven;}
   \item{else ngiven2 =nng2;}
\end{itemize}
end;\\

\subsection{e. Determination of the properties of h.} This is a rather long 
section of Autsubc. The bullets, dashes and asterisks below are an artifact 
of the way the indentation in the LaTeX program is set up and are not in the 
original CAYLEY program. Note if nng is not equal to one, we just get the list 
of normal subgroups of Aut($G$) from this part of Autsubc.\\
n=normal subgroups(h);\\
L=length(n);\\
print(L);\\
set batch=false;\\
success=false;\\
{\itshape \textquotedblleft This loop runs through the proper normal subgroups of g seeking a possible first factor."}\\
for i=2 to L-2 do
\begin{itemize}
   \item{ord=fetch(n,i;order);}
   \item{if(ord ge n4) then}
    \begin{itemize} \item{ for j=i+1 to L-1 do}
     \begin{itemize} \item{ if(order(n[i])*order(n[j]) eq order(h)) then}
        \item{if((n[i] meet n[j]) eq $<$ identity of h $>$) then}
        \item{print 'h is a direct product of:';}
        \item{print i,j;}
         \item{print classes(n[i]);}
          \item{print classes(n[j]);}
          \item{success=true;}
          \item{break;}
        \end{itemize}\item{end;}
      \end{itemize}\item{end;}
   \end{itemize} end; "second factor"
   \begin{itemize}\item{if(success) then}
       \item{break;}
     \item{end;}
     \item{end;}
    \end{itemize} end; "first factor"
\begin{itemize}\item{if(not success) then}
  \item{print 'h is not a direct product';}
\item{print classes(h);}
\item{"print n;"}
  \end{itemize} end;\\
invh = invariant(autg,h);
\begin{itemize}\item{if(invh) then}
  \item{ print invariant(autg,h);}
   \item{elh=elementary abelian(h);}
   \item{print elh;}
\end{itemize}
else
\begin{itemize}
\item{print invariant(autg,h);}
\item{print classes(h);}
\item{quit;}
\end{itemize}
end;\\
if(order(h) ne orf) then
\begin{itemize}
   \item{print order(h);}
   \item{quit;}
\end{itemize} end;\\
set batch=false;\\
nq=0;\\
np=0;\\
np2=0;\\
np3=0;\\
np4=0;\\
np5=0;\\
np6=0;\\
npq=0;\\
np2q=0;\\
np3q=0;\\
np4q=0;\\
np5q=0;\\
np6q=0;\\
lq=0;\\
lp=0;\\
lp2=0;\\
lp3=0;\\
lp4=0;\\
lp5=0;\\
lp6=0;\\
lpq=0;\\
lp2q=0;\\
lp3q=0;\\
lp4q=0;\\
lp5q=0;\\
lp6q=0;\\
set batch=true;\\
{\itshape \textquotedblleft Print classes(h).";}\\
cl=classes(h);\\
lc=length(cl);\\
print lc;\\
set batch=false;\\
for i=1 to lc do
\begin{itemize}  \item{LL = fetch(cl,i;length);}
  \item{ox = fetch(cl,i;order);}\\
  if(ox eq q) then
    \begin{itemize} \item{nq=nq+1;}
    \item{lq=lq+LL;}
    \end{itemize}
   \item{end;}\\
   if(ox eq p) then
     \begin{itemize} \item{np=np+1;}
     \item{lp=lp+LL;}
   \end{itemize} \item{end;}\\
   if(ox eq p$^{\wedge} $2) then
      \begin{itemize} \item{np2=np2+1;}
      \item{lp2=lp2+LL;}
    \end{itemize} \item{end;}\\
  if(ox eq p$^{\wedge} $3) then
      \begin{itemize} \item{np3=np3+1;}
      \item{lp3=lp3+LL;}
    \end{itemize} \item{end;}\\
  if(ox eq p$^{\wedge} $4) then
     \begin{itemize} \item{np4=np4+1;}
      \item{lp4=lp4+LL;}
    \end{itemize} \item{end;}\\
  if(ox eq p$^{\wedge} $5) then
      \begin{itemize} \item{np5=np5+1;}
      \item{lp5=lp5+LL;}
    \end{itemize} \item{end;}\\
  if(ox eq p$^{\wedge} $6) then
      \begin{itemize} \item{np6=np6+1;}
      \item{lp6=lp6+LL;}
    \end{itemize} \item{end;}\\
   if(ox eq p*q) then
     \begin{itemize}\item{npq=npq+1;}
     \item{lpq=lpq+LL;}
   \end{itemize}\item{end;}\\
   if(ox eq p$^{\wedge} $2*q) then 
      \begin{itemize}\item{np2q=np2q+1;}
      \item{lp2q=lp2q+LL;}
    \end{itemize}\item{end;}\\
  if(ox eq p$^{\wedge} $3*q) then
      \begin{itemize} \item{np3q=np3q+1;}
      \item{lp3q=lp3q+LL;}
    \end{itemize}\item{end;}\\
  if(ox eq p$^{\wedge} $4*q) then
      \begin{itemize} \item{np4q=np4q+1;}
      \item{lp4q=lp4q+LL;}
    \end{itemize}\item{end;}\\
  if(ox eq p$^{\wedge} $5*q) then
      \begin{itemize} \item{np5q=np5q+1;}
      \item{lp5q=lp5q+LL;}
    \end{itemize} \item{end;}\\
  if(ox eq p$^{\wedge} $6*q) then
      \begin{itemize} \item{np6q=np6q+1;}
      \item{lp6q=lp6q+LL;}
    \end{itemize} \item{end;}
\end{itemize}end;\\
set batch=true;\\
z=center(h);\\
print z,order(z);\\
zr=relations(z);\\
print zr;\\
qq=h/z;\\
print nclasses(qq);\\
{\itshape \textquotedblleft order structure of h"}\\
set batch=false;\\
{\itshape \textquotedblleft order of element, number of elements, number of classes"}\\
print q, lq, nq;\\
print p, lp, np;\\
print p$^{\wedge} $2, lp2, np2;\\
print p$^{\wedge} $3, lp3, np3;\\
print p$^{\wedge} $4, lp4, np4;\\
print p$^{\wedge} $5, lp5, np5;\\
print p$^{\wedge} $6, lp6, np6;\\
print p*q, lpq, npq;\\
print p$^{\wedge} $2*q, lp2q, np2q;\\
print p$^{\wedge} $3*q, lp3q, np3q;\\
print p$^{\wedge} $4*q, lp4q, np4q;\\
print p$^{\wedge} $5*q, lp5q, np5q;\\
print p$^{\wedge} $6*q, lp6q, np6q;\\
set batch=true;\\
\textquotedblleft u=automorphism group(h);\\
print order(u),ngenerators(u),degree(u);"\\
\subsection{f. Properties of Aut($G$)} The following part deals with mappings 
of the normal subgroups of Aut($G$). We need the automorphism group of 
Aut($G$) here. This is the group v below. Here we are still dealing with the 
assumption that h is the subgroup of Aut($G$) that we are interested in. This 
can be seen from the line [x eq h] in the program below.\\

\noindent v=automorphism group(autg);\\
print 'properties of aut(aut(g))';\\
print forder(v),order(v),ngenerators(v),degree(v);\\
ng=ngenerators(v);\\
m=normal subgroups(autg);\\
set batch=false;\\
L=length(m);\\
print 'length of m =', L;\\
for i=2 to L-1 do
\begin{itemize}
   \item{ord = fetch(m,i;order);}\\
   if(ord eq ngiven) then\\
     \begin{itemize} \item{x=m[i];}
     \item{ncl=nclasses(x);}
     \item{ z=center(x);}
     \item{ r=relations(x);}
     \item{ print '-------------------';}
     \item{ print i,order(m[i]);}
     \item{ print 'number of classes=',ncl,'center relations', r;}
     \item{XZ = x eq h;}
     \item{ print XZ;}
     \end{itemize}
     if(ncl eq lc) then
      \begin{itemize}
        \item{print relations(x);}
     \end{itemize}
\end{itemize}
nn=normal subgroups(x);\\
Ln=length(nn);\\
set batch=false;\\
success=false;\\
{\itshape \textquotedblleft This loop runs through the proper normal subgroups of g seeking a possible first factor."}\\
for i1=2 to Ln-2 do
   \begin{itemize} \item{ord=fetch(nn,i1;order);}
  \item{ if(ord ge n4) then}
     \begin{itemize} \item{for j1=i1+1 to Ln-1 do}
      \begin{itemize}\item{if(order(nn[i1])*order(nn[j1]) eq order(x)) then}
       \item{if ((nn[i1] meet nn[j1]) eq  $<$ identity of x $>$) then}
        \item{print 'x is a direct product of:';}
          \item{print i1,j1;}
          \item{print classes(nn[i1]);}
          \item{print classes(nn[j1]);}
          \item{success=true;}
          \item{break;}
        \end{itemize}\item{end;}
      \end{itemize} \item{end;}
     \end{itemize} end; "second factor" \\
   if(success) then
       \begin{itemize} \item{break;}
    \end{itemize} end; \\
     end;\\
     end; "first factor"\\
if(not success) then
  \begin{itemize} \item{print 'x is not a direct product';}
\item{if(XZ) then}
  \begin{itemize}\item{print classes(x);}
  \end{itemize} \item{end;}
  \item{end;}
\end{itemize} 
end;\\
end;\\
end;\\
"s2=sylow subgroup(autg,2);\\
XZ=h eq s2;\\
print XZ;"\\
\subsection{g. Mappings of the normal subgroups of Aut($G$)} This part of the 
program is most useful if our group h is not the group of interest, i.e., h 
is not the group associated with the extension.\\
chars=empty;\\
print 'the length of the normal subgroups of the aut group factors is ', l;\\
print 'the order of the aut(aut(g) factor) is ', order(autg);\\
nautgens=ngenerators(v);\\
{\itshape \textquotedblleft Determine which normal subgroups of G are left invariant by the automorphism group."}\\
print 'the characteristic subgroups of aut factor are';\\
for i=2 to l-1 do\\
\begin{itemize}\item{   invar = true;}
   \item{for j=1 to ng do}
     \begin{itemize} \item{f=automorphism(aut,v.j);}
     \item{invar = f(m[i]) eq m[i];}
     \item{if(not invar) }
        \begin{itemize} \item{then break;}
         \end{itemize} \item{end;}
         \end{itemize} \item{end; "automorphism group generators"}
         \item{if(invar) then}
         \begin{itemize} \item{print i,order(m[i]);}
          \item{chars = append(chars,m[i]);}
           \end{itemize} 
          \item{end};
          \end{itemize} 
    end; "normal subgroups"\\
print 'end of characteristic subgroups';\\
print '-----------------------------';\\
for i=2 to l-1 do
\begin{itemize}
   \item{nequal=m[i] eq h;}
   \item{if(nequal) then}
      \begin{itemize} \item{print 'h is equal to m[',i,']'; }
   \end{itemize} \item{end;}
  \end{itemize} end;\\
print 'mappings of noncharacteristic subgroups of aut 2 group';\\
print 'among themselves by aut(aut(2 group))';\\
print 'the number of generators in aut(aut(32 group))=',ng;\\
for i=2 to l-1 do
\begin{itemize}
   \item{ord = fetch(m,i;order);}
   \item{if(ord eq ngiven2) then }
      \begin{itemize} \item{invar = false;}
      \item{for j=1 to ng do}
         \begin{itemize} \item{f=automorphism(autg,v.j);}
         \item{X=f(m[i]);}
         \item{invar = X eq h;}
         \item{if(invar) then}
           \begin{itemize} \item{print 'for i=',i,'f(m[i]) under 
autg(v.',j,') equals h ' ;}
           \item{else}
            \item{for k=2 to l-1 do}
               \item{if(X eq m[k]) then;}
                \item{if(i ne k) then }
                \item{orm = order(m[k]);}
            \item{print 'group m[',i,'] mapped by v[',j,'] to m[',k,'] of 
order', orm;}
            \item{end;  "search for image of m under f" }
            \item{end; "check to see if m is characteristic for u.i" }
        \item{end; "end of search under aut grp generators"}
  \end{itemize} \item{end; "end of loop on normal subgroups"}
\item{invar=false;}
\end{itemize} \item{end;}
\end{itemize} \item{end;}
\end{itemize} end;\\
set columns = 70;\\
lnc=length(chars);\\
print 'number of characteristic subgroups of aut(32 group) =',lnc;\\
for i=2 to l-1 do
\begin{itemize}
  \item{ nequal=m[i] eq h;}
   \item{if(nequal) then}
     \begin{itemize}
     \item{print 'h is equal to m[',i,']';}
  \end{itemize} \item{end;}
\end{itemize}end;\\
quit;\\
\$exit\\
\\

\section{Tables}

\begin{tabular}{|c|c|c|} \hline
  \multicolumn{3}{|c|}{Table 1a} \\ \hline
  \multicolumn{3}{|c|}{Groups of order 32 and their automorphism groups}  \\
\hline
\multicolumn{3}{|c|}{I. Abelian Cases (numbers 1 to 7)} \\ \hline
	2-group        &      presentation      &  automorphism group   \\ \hline
1.   $  (1^5)          $&$   a^2=b^2=c^2=d^2=e^2= $&      GL(5,2)         \\
		         &$   (a,b)=(a,c)=(a,d)=   $&                       \\
		        &$    (a,e)=(b,c)=(b,d)=  $&                       \\
		        &$    (b,e)=(c,d)=(c,e)=  $&                      \\
		        &$    (d,e)=1            $&                       \\ \hline
2. $ (2,1^3)     $     &$  a^4=b^2=c^2=d^2=(a,b)=$&     order 21,504     \\
      	        &$    (a,c)=(a,d)=(b,c)=  $&                      \\
      	        &$    (b,d)=(c,d)=1      $&                    \\ \hline
3. $ (2^2,1)  $       &$  a^4=b^4=c^2=(a,b)=    $ &    order 1536       \\
        	        &$    (a,c)=(b,c)=1      $&                   \\ \hline
4. $ (3,1^2)         $&$  a^8=b^2=c^2=(a,b)=    $& Aut$(C_4 \times  C_2 
\times
C_2) \times  C_2 $\\
       	        &$    (a,c)=(b,c)=1      $&                      \\ \hline
5.    (3,2)          &$  a^8=b^4=(a,b)=1      $&order 128, number 753\\
\hline
6.    (4,1)          &$  a^{16}=b^2=(a,b)=1     $&  order 32, number 10 \\
\hline
7.    (5)	        &$  a^{32}=1               $&$      C_8  \times  C_2  $  
\\
\hline

\end{tabular}

\begin{tabular}{|c|} \hline
Notes for Table 1a \\ \hline
       \#2. order 21,504.   \\ \hline
            CAYLEY returns this automorphism group as a permutation  \\
   	      group of degree 31 on 9 generators. The generators are       \\
            all order two, and the presentation given by CAYLEY can      \\
            be reduced to one with fewer generators. A five generator            \\
            presentation of this automorphism group is: \\                                     

     		$a^2=b^2=c^2=d^2=e^2=(a*d)^2=(a*e)^2=(b*e)^2=$\\
\hline
\end{tabular}

\begin{tabular}{|l|c|c|} \hline
  \multicolumn{3}{|c|}{Table 1b} \\ \hline
  \multicolumn{3}{|c|}{Groups of order 32 and their automorphism groups}  \\
\hline
\multicolumn{3}{|c|}{II. Isoclinic Class 2. (numbers 8 to 22)} \\ \hline
          2-group        &     presentation          & automorphism group
\\ \hline
8. [2,a1]$ D_4 \times  (1^2)$ & $\left[ \ast \right ]a^4=b^2=a^b*a=c^2=d^2=$
& Hol$(C_4 \times   C_2 \times  C_2)$ \\
generalized      &     $ (c,d)=(a,c)=(a,d)= $    &                  \\
   dihedral      &    $ (b,c)=(b,d)=1   $  &               \\ \hline
9. [2,a2]$ Q_2 \times  (1^2)$ & [$\ast$]$a^4=b^4=a^2*b^2=a^b*a= $&    order
9216       \\
\hline
\end{tabular}

\begin{tabular}{|c|} \hline
Notes for Table 1b. Aut($G$) Presentations.      \\ \hline
\# 9. $ Q_2 \times  C_2  \times  C_2 $'s automorphism group          \\
\hline
           $    a^2=b^2=c^2=d^4=e^2=f^4=(a,d)=(a,e)=(a,f)=(b,d)= $\\
       $	   (b,f)=(c,d)=(c,f)=(d,e)=(a*b)^3=c*d^2*e*c*e=          $\\
       $	   d^2*f*(d^{-2})*f=d*f*d*((f*d*f)^{-1})=e*f^2*e*(f^{-2})=  $\\
       $	   a*b*c*b*c*a*c=(b*e)^4=(e*f)^2*(e*(f^{-1}))^2=        $\\
       $	   a*b*a*e*b*e*a*b*e=d*e*f*d*e*(f^{-1})*e*(d^{-1})*(f^{-1})=1 $ \\
\hline
\#10. order 768.                                                      \\
\hline
\end{tabular}   \\

\begin{tabular}{|c|} \hline
Notes for Table 1b, continued. Aut($G$) Presentations.  \\ \hline

\#16. order 256. Cayley returns aut(g) as a group of degree 28 \\ 
on 8 generators. One presentation for the aut(g) is:  \\
                                                         \\
       	 $   a^2=b^2=c^2=e^2=f^2=g^2=h^2=a*d^2= $ \\
       	 $   (a,b)=(a,c)=(a,e)=(a,f)=(a,g)=(a,h)= $\\
       	 $   (b,c)=(b*d)^2=(b,e)=(b,f)=(b,g)=(b,h)= $\\
       	 $   (c*d)^2=(c,e)=(c,g)=(c,h)=(d,f)=(d*(h^{-1}))^2=  $\\
\hline
\end{tabular}\\
\newpage
\begin{tabular}{|l|c|c|} \hline
  \multicolumn{3}{|c|}{Table 1c} \\ \hline
  \multicolumn{3}{|c|}{Groups of order 32 and their automorphism groups}  \\
\hline
\multicolumn{3}{|c|}{III. Isoclinic Class 3. (numbers 23 to 32)} \\ \hline
      2-group        &     presentation          & automorphism group   \\
\hline
23. $[3,a1] D_8 \times  C_2 $&  $ a^8=b^2=c^2=(a,b)*a^2= $&  Hol($C_8 \times
  C_2$)  \\
                       & $   (a,c)=(b,c)=1 $      &                     \\
       (generalized dihedral)    &                          &                      \\
\hline
24. $[3,a2]$          &$  b^2=c^2=(b,a)*a^2=  $&   \#2216 order 128 \\
\hline
\end{tabular}\\
\\

\begin{tabular}{|l|} \hline
\multicolumn{1}{|c|}{Notes for Table 1c.} \\ \hline
 \#23. Hol($C_8 \times  C_2)$ = $(C_8 \times  C_2) @ (D_4 \times
C_2)$: \\
   \qquad \qquad $ a^8=b^2=(a,b)=c^4=d^2=c^d*c=e^2=(c,e)=(d,e)=  $
        \\
    \qquad \qquad  $ a^c*b*(a^{-1})=b^c*b*(a^{-4})= a^d*a^3 =  $
      \\
   \qquad \qquad      $ b^d*b*(a^{-4})= a^e*a = (b,e) = 1  $
    \\ \hline
\#24. The presentation returned by CAYLEY is: \\
\\
       	\qquad \qquad   $ a^2=b^2=c^2=d^4=e^2=(a,c)=(a*d)^2=(a,e)=(b,c)= $            \\
\hline

\end{tabular}\\
\\

\begin{tabular}{|l|c|c|} \hline
  \multicolumn{3}{|c|}{Table 1d} \\ \hline
  \multicolumn{3}{|c|}{Groups of order 32 and their automorphism groups}  \\
\hline
\multicolumn{3}{|c|}{IV. Isoclinic Class 4. (numbers 33 to 41)} \\ \hline
      2-group        &     presentation          & automorphism group   \\
\hline
33. [4,a1]          &$  a^2=b^2=c^2=(a,b)=      $&     order 384 \#20,100     
   \\
      $ (1^2) \wr C_2$ &$    ((a,c),a)=((b,c),a)=  $&  $ C_2 \wr S_3$
  \\
                      &$    ((b,c),b)=1         $&  \\
\hline
34. [4,a2]         &$  a^4=b^4=c^2=(a,b)=     $&      Hol($C_4 \times  C_4$)
     \\
        Dih($2^2$)    &  $  (a,c)*a^2=(b,c)*b^2=1 $&                   \\
\hline
35. [4,a3]         &$  b^4=a^2*(c^{-2})=(a,b)=  $& order 512, not wreath \\
                       &$  (c,a)*a^2=(b,c)*b^2=1 $&  product \\
\hline
36. [4,b1]  &$ a^4=b^2=c^2=(a,c)*a^2= $& $ (1^2) \wr C_2 \times  (1^2) $  \\
      $ <2,2|4;2> @ C_2$& $     (a,b)=((b,c),c)=1    $&                 \\
\hline

37. [4,b2] &$ a^2*(c^{-2})=b^2=(c,a)*a^2= $& $ (1^2) \wr C_2 \times (1^2)$\\
                      &$      (a,b)=((b,c),b)=1   $&
\\ \hline

   38. [4,c1] &$ a^2=c^2=(c,a)*b^2=(a,b)=$&$ (C_2 \times  C_2 \times C_2) \wr C_2$ \\
                      &$      ((b,c),c)=1         $&  \\
\hline
39. [4,c2]         &$ b^4=c^2=(c,a)*a^2=(a,b)=  $&       order 256 \\
      $(C_4 \times  C_4) @ C_2  $&$    (b,c)*b^2*(a^-2)=1 $& \#16,888              
\\
\hline
40. [4,c3]         &$   a^2*(c^{-2})*(b^{-2})=   $&      order 256     \\
        	           &$   (c,a)*a^2=(b,c)*c^2=   $&   \#5675                  
   \\
    	                 &$        (a,b)=1          $&                   \\
\hline
41. [4,d]          &$  c^2=(c,b)*a^2=     	 $& $  1^6 @ C_3$ same as     \\
        	           &$   (c,a)*a^2*b^2=         $&$  (1^2) \wr C_3$, and
   \\
                      &$      (a,b)=1  $&$ (1^2) \times (1^4) @ C_3$  \\
\hline

\end{tabular}\\
\\

\begin{tabular}{|l|c|c|} \hline
  \multicolumn{3}{|c|}{Table 1d continued} \\ \hline
  \multicolumn{3}{|c|}{Groups of order 32 and their automorphism groups}  \\
\hline
\multicolumn{3}{|c|}{V. Isoclinic Class 5. (numbers 42 to 43)} \\ \hline
      2-group        &     presentation          & automorphism group   \\
\hline
42. [5,a1]	    &$ a^2*(b^{-2})=a^2*(c^{-2})=    $&  $ S_4 \wr C_2 $\\
    $Q_2 \text{Y} Q_2 = D_4 \text{Y} D_4$ &$ a^2*(d^{-2})=(c,b)*a^2= $&      
           \\
                      &$(d,a)*a^2=(a,b)=(a,c)=    $&
\\
   (extra-special)    &$(b,d)=(c,d)=1            $&                    \\
\hline

   43. [5,a2]	  &$ d^2=a^2*(b^{-2})=a^2*(c^{-2})=$&        order 1920    \\
      $ Q_2 \text{Y} D_4$    &$ (c,b)*a^2=(d,a)*a^2=(a,b)=$&      complete 
group
\\
                      &$(a,c)=(b,d)=(c,d)=1 $      & \#240,996
\\
     (extra-special)  &                              &                    \\
   Dirac Matrix group &                               &
\\ \hline

\end{tabular}\\
\\
\newpage
\begin{tabular}{|l|} \hline
\qquad \qquad \qquad \qquad \qquad Notes for Table 1d. \\ \hline
    \#33. We do not have a simple description in terms of permutations. \\
       	\qquad A presentation is
         \\
       	  $ a^2=b^2=c^2=d^3=(a*b)^2=(c*d)^2=a*c*(d^{-1})*a*d*c=$ \\
       	$ (a*d)^3=(a*c*b*c)^2=(b*c)^4=(b*d*b*(d^{-1}))^2= $\\                            
       	$ b*c*b*d*b*d*c*b*(d^{-1})=1$ [same as  \# 18 ?]. \\ \hline
\#35. order 512. This automorphism group has 53 classes with         \\
       	 \qquad the order structure:                           \\
\hline

\end{tabular}\\
\\

\begin{tabular}{|l|c|c|} \hline
  \multicolumn{3}{|c|}{Table 1e} \\ \hline
  \multicolumn{3}{|c|}{Groups of order 32 and their automorphism groups}  \\
\hline
\multicolumn{3}{|c|}{VI. Isoclinic Class 6. (numbers 44 to 45)} \\ \hline
      2-group        &     presentation          & automorphism group   \\
\hline
  44. [6,a1]        &$  a^2=c^2=(b,a)*b^4= $&$  D_4 \times D_4 $ \\
        Hol($C_8$)      &$ (b,c)*b^2=(a,c)=1 $& \#154 order 64 \\ \hline
   45. [6,a2]        &$   a^2=b^4*(c^{-2})=(b,c)*b^2= $&$   D_4 \times  D_4$ \\
                     &$    (b,a)*b^4=(a,c)=1     $&     \#154 order 64   \\
\hline
\end{tabular}\\
\\

   \begin{tabular}{|c|} \hline
Notes for Table 1e. \\ \hline
    The relations for \#49, \#50, and \#51 are different from those given by Sag \\
 and Wamsley. This does not mean that their relations in these cases are   \\
     wrong, but just that these are simpler ones than theirs are to deal with.
\\ \hline
               Presentations for automorphism groups                     \\
\hline
    \#46. Hol($C_4 \times  C_2$) \# 259 of order 64. Relations from Sag and
Wamsley. \\
       	$a^2=b^2=c^2=(a,b)^2=(a,c)^2=(b,c)=((a,b),c)*(((a,c),b)^{-1})=1$ \\
\hline
  \#47. order 128.                                                      \\
       	   $ a^2=b^2=c^4=(a,b)=(b*c)^2=(a*c)^4=$ \\  
       		$a*c^2*a*c^2*a*c^{-2}*a*c^{-2}=1 $ \\                      \\
       	    In this presentation the group $h=<b,c>$ is just             \\
       	    $D_4$. This automorphism group is $D_4 \wr C_2$.   \\ \hline
    \#48. order 128.                                                    \\ 
              $ a^2=b^2=c^4=(a,b)=b*c^2*b*(c^{-2})=(a*c)^4=  $ \\
       		$(a*c*b*c)^2=(b*c)^2*(b*(c^{-1}))^2=                    $ \\
		$(a*b*c)^2*a*(c^{-1})*a*c=1                               $\\
       	One of the normal subgroups of this order 128 group            \\
  	        	is  $D_4 \times  D_4$. [ Aut is isomorphic to $D_4 \wr C_2.$ ] \\

\\ \hline
       \#50 Number 117 of order 64. \\
      	      $  a^4=b^8=c^2=(a,b)*b^4=(b,c)*b^2=(a,c)=1          $\\
\hline
\end{tabular}\\

\begin{tabular}{|c|c|c|c|c|c|c|c|} \hline
  \multicolumn{8}{|c|}{Table 2a} \\ \hline
  \multicolumn{8}{|c|}{Groups of order 32$p$ and their automorphism groups}  \\
\hline
\multicolumn{8}{|c|}{I. Abelian groups of order 32} \\
\multicolumn{8}{|c|}{The number of groups of order 32 = 7.} \\
\multicolumn{8}{|c|}{The number of groups of order $32p$ with $C_2$ images = 12.}
\\ \hline
  &  $C_2$ & \multicolumn{4}{|c|}{image of 2-group} & normal  & group \\
\cline{3-6}
2-group & action& $C_2$ &$C_4$ & $C_8$& $C_{16}$ &2-group & associated \\
  & & & & & & for  & with  \\ 
  & & & & & &   $p=3$ $\dagger $& extension \\ \hline
  \#1. &  $a$ 31 & Hol($1^4$)  & &   &  &$ (1^3)\times  A_4$ &
(1,1,1,1)\\
$1^5$&    &    &     &   &  & $C_2 \times  $ [48]  &  \\ \hline
   \#2.&  $a$ * & Aut$(2,1^3)$ & $a$ Hol($1^3$)& &  & (2,1) $\times  A_4$
& (1,1,1,1)\\
  $(2,1^3)$&  $b$ 14 &order 1536 &          &    &    &            &
(2,1,1) \\ \hline
  \#3. & $a$ 3  &order 512  & $a$ Hol(2,1)&  &  & [48]$ \times  C_2 $&
(2,1,1) \\
  $(2^2,1)$&  $c$ 4  &order 384  &           &   &    &            &$
(2^2)$\\ \hline
  \#4.&  $a$ *  &Aut(211)$\times C_2$& $a$   [192] &  $S_4$ &  & $ C_8 
\times
A_4 $&  (2,1,1) \\
$(3,1^2)$ & $b$ 6  &  \#33$ \times  C_2$ &$a,b$ Hol(2,1)&     &    &         
&
(3,1)  \\ \hline
  \#5. &  $a$ *  &  Aut(3,2) & $a$   \#33   &$  D_4$ &   &    &  $(2^2)$
  \\
   (3,2)  &  $b$ 2  & \#33 $\times  C_2$& $b$  $D_4 \times  C_2$&   &    &   
  &
(3,1) \\ \hline
   \#6. &  $a$ *  &$C_4\text{Y}Q_2 \times  C_2 $& $a$ $ C_4 \text{Y} Q_2$& $ D_4$ 
&$ C_2$  &
  &  (3,1)  \\
    (4,1) &  $b$ 2  &  (2,1,1)  &$a,b$ $C_4 \text{Y} Q_2$&$  D_4$ &  &     &
  (4) \\ \hline

    \#7.   &  $a$ *  &$  C_8 \times  C_2 $& $a$  $  C_8$ &$  C_4 $&$
C_2$ &     &   (4)    \\ 
         $C_{32}$   &         &   &            &     &  &   & \\ \hline
\end{tabular}\\

\begin{tabular}{|c|} \hline
Notes for Table 2a. \\ \hline
If $p\equiv 1$ mod(32), e.g., $p=97$. then Aut($C_{97}$) = Hol($C_{97}$), \ldots . \\
\hline
An * in the $C_2$-action column means the group associated with the \\
extension is a characteristic subgroup. \\

A number in the $C_2$ column gives the multiplicity of the group associated  \\
with the extension. \\

In column 7, we give the groups of order 96 with the 2-subgroup in  \\
column 1 as its normal subgroup. \\
\hline
\end{tabular}\\
\\

\begin{tabular}{|c|c|c|c|c|c|} \hline
  \multicolumn{6}{|c|}{Table 2b} \\ \hline
  \multicolumn{6}{|c|}{Groups of order 32$p$ and their automorphism groups}  \\
\hline
\multicolumn{6}{|c|}{Isoclinic Class 2  } \\
\multicolumn{6}{|c|}{The number of groups of order 32 = 15.} \\
\multicolumn{6}{|c|}{The number of groups of order 32$p$ with $C_2$ images = 42.}
\\ \hline
  &  $C_2$ & \multicolumn{2}{|c|}{image of 2-group} & normal  & group \\
\cline{3-4}
2-group & action& $C_2$ &$C_4$ & 2-group & associated \\
  & & && for $p=3\dagger $ &  with  \\
\hline
\end{tabular}

\begin{tabular}{|c|c|c|c|c|c|c|} \hline
  \multicolumn{7}{|c|}{Table 2b continued} \\ \hline
  \multicolumn{7}{|c|}{Groups of order 32$p$ and their automorphism groups}  \\
\hline
\multicolumn{7}{|c|}{Isoclinic Class 2  } \\ \hline
  &  $C_2$ & \multicolumn{3}{|c|}{image of 2-group} & normal  & group \\
\cline{3-5}
2-group & action& $C_2$ &$C_4$ & $C_8$&  2-group & associated \\
  & & & & &  for $p=3\dagger $ & with  \\ 
  & & & & &     & extension \\ \hline
   16.    & $c$  *  & Aut(16)   & $c$  $Q_2\text{Y}Q_2 $ & &  & 
(2,1,1)  \\
	     & $a$  2  & $1^3 \wr C_2$ &        &  & &   
$(4,4|2,2)$\\
          & $b$  2  & $1^3 \wr C_2$ &        & &  &   
$C_4 \times  C_4 $ \\
          &$a,b$ 2  & \#33 $\times  1^2 $&        &  &       & $C_4
@ C_4 $ \\ \hline
  17.    & $a$  *  & $S_4 \times  1^2 $ &$a,b,c$ $S_4 \times  C_2$&  & 
[48]$ \times  C_2 $ & $C_4 \text{Y} Q_2 $ \\
  $C_8 \text{Y} Q_2$ & $c$  3  & $D_4 \times  1^2 $ & $a$ $ D_4 \times  
C_2$& 
&            & $C_8 \times  C_2 $ \\
          &$a,b$ 3  & $D_4 \times  1^2 $ &        &  &  &         
$<2,2|2> $ \\ \hline
18.    & $a$  3  & $1^3 \wr C_2$ & $b$  \#33    &  &     1 &  
(2,1,1)  \\ \hline
   19.    & $b$  *  &  Aut(19)  & $b$  \#33    &  &             & $C_4
\times  C_4 $ \\
          & $a$  2  & \#33 $\times  C_2$  & $a$  $D_4 \times  C_2$&  & 
          & $C_8 \times  C_2 $ \\ \hline
   20.    & $b$  *  & \#33$ \times  C_2$  & $b$  \#33    &$b$ $D_4$ & 
    & (2,1,1)  \\
          & $a$  2  &   $1^5$     & $a,b$ \#33    &  &             & $C_8
\times  C_2 $ \\ \hline
   21.    & $b$  *  & \#33$ \times  C_2$  & $b$  \#33    &$b$ $D_4$ & 
    & $C_4 \times  C_4$  \\
$(4,8|8,8)$& $a$  2  &   $1^5 $    & $a,b$ \#33    &  &  &  $C_8
\times  C_2$  \\ \hline

   22.    & $b$  *  &\#8[16]$ \times  C_2$& $b$  $C_4 \text{Y} Q_2$&$b$  
$D_4$& & $C_8 \times  C_2 $ \\
          & $a$  2  & (2,1,1)   & $a,b$ $C_4 \text{Y} Q_2$&$ab$ $D_4$&  & $C_{16}$  \\ \hline

\end{tabular}\\
\\

\begin{tabular}{|c|} \hline
Notes to Table 2b. \\ \hline
See the general comments in the notes to Table 2a. They apply to all of\\
          Table 2 (2a, 2b, \ldots to 2e).                                  \\
\hline
\#8a. Relations for the automorphism group of order 1536:           \\
           $     a^2=b^2=c^2=d^2=e^2=(a,b)=(c,d)=(a*c)^3= $
\\
            $   (a*d)^3=(a*e)^3=(b*c)^3=(b*d)^3=           $
\\
             $  (b*e)^3=(c*e)^3=(d*e)^3=(a*c*a*e)^2=        $
\\
              $ (a*d*a*e)^2=(a*e*b*c)^2=(b*d*b*e)^2=1       $
\\
\hline
\end{tabular}\\

\begin{tabular}{|c|} \hline
Notes to Table 2b continued \\ \hline

        \#10b. Same as \#10a.
\\ \hline
        \#10a,b,c factor is order 768 \#1087581 \\ \hline
        \#11a. The relations for the order 256 (\# 8935) group are:
\\

\\
          $      a^2=b^2=c^2=d^2=(a*d)^2=(b*c)^2=(a*d)^4= $
\\
           $     (a*b*d*b)^2=(a*b*a*c)^2=(a*c)^4=          $
\\
            $    (a*c*d*c)^2=(b*d)^4=(b*d*c*d)^2=(c*d)^4=1 $
\\

\\
           This automorphism group has 247 normal subgroups and 58        \\
            conjugacy classes. The class structure is:
\\
               order 2 [$1^{15},4^{24}$]
     \\
               order 4 [$8^{18}$]
   \\

\\
               Representation in terms of Aut(32 \#11) generators:
\\
               $<b^2,(b*d)^2,d,a*b*c*b*a*b*c*d*b,(a*b)^2,$
\\
                       $(b*c)^2,a,c>$; order 256.
\\ \hline
        \#12b. Same as \#11a. This is \#8935 of order 256 in Small Group 
Library,
\\
                Representation in terms of Aut(32 \#12) generators:
\\
                $<c^2,(c*d)^2,d,(a*e)^2,a,b,(c*e)^2,e>$;
\\
                        order 256.
\\ \hline
        \#15c. The order 96 group is the one in Table 3 IIC \#36.
\\

\\
               $(1^4 @ C_3) @ C_2$. The automorphism group of this group
   \\
                                has order 576.
\\ \hline
        \#16a. Generators as a subgroup of Aut(32 \#16):
\\
                $<a^2,(a*b*d)^2,b*d,(a*b)^2,(b*c)^2,b,$
\\
                       $ b*c*b*a*c>$;  order 128.
\\ \hline
        \#16b. Same as \#16a:  $ (C_2 \times  C_2 \times  C_2) \wr C_2$.
                   \\ \hline
\end{tabular}\\
\\

\begin{tabular}{|c|c|c|c|c|c|} \hline
  \multicolumn{6}{|c|}{Table 2c} \\ \hline
  \multicolumn{6}{|c|}{Groups of order 32$p$ and their automorphism groups}  \\
\hline
\multicolumn{6}{|c|}{Isoclinic Class 3  } \\
\multicolumn{6}{|c|}{The number of groups of order 32 = 10.} \\
\multicolumn{6}{|c|}{The number of groups of order 32$p$ with $C_2$ images = 30.}
\\ \hline
  &  $C_2$ & \multicolumn{2}{|c|}{image of 2-group} & normal & group \\
\cline{3-4}
2-group & action& $C_2$ &$C_4$ &   2-group & associated \\
  & & & &  for $p=3\dagger $ & with  \\
\hline

\end{tabular} \\

  \begin{tabular}{|c|} \hline
Notes for Table 2c. \\ \hline
Many of the general comments made in the notes to Table 2a also  \\
       	 apply to this table.  \\
\hline
    \#23a. Relations for this group of order 128 \#2216 are:            \\
       	   	     $ a^2=b^2=c^2=d^2=e^2=(a,b)=(a,c)= $                   \\
       		     $ (a,d)=(b,c)=(b,d)=(c,e)=(d,e)=   $                   \\
       		     $ (a*e)^4=(c*d)^2*(e*b)^2=1;       $                   \\
       	    This automorphism group has 432 normal subgroups and 44     \\
       	    conjugacy classes. The order structure of this group is     \\
       			        order 2 [$1^3,2^{14},4^8$] \\
       			        order 4 [$2^4, 4^{14}$]  \\
       		      This group is a characteristic subgroup of \\
       		      Hol$(C_8 \times  C_2)$. This is isomorphic to \\
       		      Aut(\#24 of order 32). In terms of the generators    \\
       		      of Aut(32 \#24) this group can be represented as:   \\
        		      $<b,c*b*(c^{-1}),c^2,d,d*c,a^2*b*d,e>$ \\ \hline
         	    \#24c. The group $D_4 \times  C_2 \times  C_2$ appears many
times as          \\
       		      a normal subgroup of Aut(32 \#24). The structure     \\
       		      is as follows:  \\
       		         twice as a characteristic subgroup \\
      		         twice with multiplicity 8  \\
                          twice with multiplicity 4 \\
                          once with multiplicity 2. \\
       		      The one with multiplicity four is the one that      \\
       	         arises as the automorphism group factor here.  \\ \hline
       	    \#25a. same as \#23a (i.e., order 128 \#2216).                   
                       \\
\hline

\end{tabular}\\
\\

\begin{tabular}{|c|c|c|c|c|c|c|c|} \hline
  \multicolumn{8}{|c|}{Table 2d} \\ \hline
  \multicolumn{8}{|c|}{Groups of order 32$p$ and their automorphism groups}  \\
\hline
\multicolumn{8}{|c|}{Isoclinic Class 4  } \\
\multicolumn{8}{|c|}{The number of groups of order 32 = 9.} \\
\multicolumn{8}{|c|}{The number of groups of order 32$p$ with $C_2$ images = 33.}
\\ \hline
  &  $C_2$ & \multicolumn{4}{|c|}{image of 2-group} & normal & group \\
\cline{3-6}
2-group & action& $C_2$ &$C_4$ & $C_8$& $C_{16}$ & 2-group & associated \\
  & & & & & & for $p=3 \dagger$ & with extension \\
\hline
\end{tabular}\\

\begin{tabular}{|c|} \hline
Notes for Table 2d. \\ \hline
\#33ac. Same as \#33a. \\ \hline
\#34b. The relations here are the same as \#35b  (256 \#16,888).   \\ \hline
  	\#35b. In terms of the generators of Aut(32 \#35) we have         \\
        		$<a^2,b*a^2*(b^{-1}),b^2,(a*c)^2,d,c*d,a*b*c*a*c*b>$          \\
       			order 256 (\#16,888). \\
       	      relations: \\
       		   $  a^2=b^2=c^2=d^2=(a*b)^2=(a*c)^2= $  \\
\hline

\end{tabular}\\
\\

\begin{tabular}{|c|c|c|c|c|c|c|} \hline
  \multicolumn{7}{|c|}{Table 2e} \\ \hline
  \multicolumn{7}{|c|}{Groups of order 32$p$ and their automorphism groups}  \\
\hline
\multicolumn{7}{|c|}{Isoclinic Class 6  } \\
\multicolumn{7}{|c|}{The number of groups of order 32 = 2.} \\
\multicolumn{7}{|c|}{The number of groups of order 32$p$ with $C_2$ images = 10.}
\\ \hline
  &  $C_2$ & \multicolumn{3}{|c|}{image of 2-group} & normal & group \\
\cline{3-5}
2-group & action& $C_2$ &$C_4$ & $C_8$& 2-group & associated \\
  & & &  & & for $p=3\dagger$ & with extension \\
\hline

\end{tabular}

\begin{tabular}{|l|c|} \hline
  \multicolumn{2}{|c|}{Table 3a} \\ \hline
  \multicolumn{2}{|c|}{Groups of order $32p$ without a normal sylow
$p$-subgroup}  \\ \hline
\multicolumn{2}{|c|}{I. Cases with a normal sylow 2-subgroup (order 96
cases)}\\  \hline
\multicolumn{2}{|c|}{A. Direct products from order 8 groups and their
automorphism groups}\\ \hline
Group & Automorphism Group \\ \hline

1.$ A_4\times C_2\times C_2\times C_2$            &   $  S_4\times
\text{Hol}(C_2\times C_2\times C_2)$     \\
                                  &                                      \\
        2.$ A_4\times C_4\times C_2$        &        $   S_4\times D_4     $
              \\
\hline
\end{tabular}\\

\begin{tabular}{|l|} \hline
\multicolumn{1}{|c|}{  Notes for Table 3a} \\ \hline
\multicolumn{1}{|c|}{Relations for the non-direct-product cases, and some 
comments}     \\
\multicolumn{1}{|c|}{on distinguishing these groups from one another.        
    }       \\ \hline
\multicolumn{1}{|c|}{IC. Cases with a normal sylow 2-subgroup. Note only the 
$C_3$ }     \\
\multicolumn{1}{|c|}{ actions are given on the 2-group with defining    }    
      \\
\multicolumn{1}{|c|}{ relations in Table 2.         }                        
      \\ \hline
         \# 11. $e^3=a^e*(b^{-1})*(a^{-1})=b^e*a=c^e*d=d^e*c*d=1 $          
   \\
                                                                         \\
\hline
\end{tabular}
\linebreak

\begin{tabular}{|c|} \hline
Notes for Table 3a continued \\\hline
Relations for Aut(\#11), Aut(\#13) and Aut(\#18) \\ \hline
               Relations for the automorphism group of \#11:               \\
                      $ a^3=b^4=c^3=d^4=(a*b)^2=(b*(d^{-1}))^2= $            
\\
                      $ a*b*c*(b^{-1})*(a^{-1})*c=       $                   
   \\
                      $ 
a*b*(d^{-1})*a*d*(b^{-1})=a*(b^{-1})*c*b*(a^{-1})*c=$      \\
                      $ a*(c^{-1})*(d^{-2})*(a^{-1})*c=                      
$   \\
                      $ 
(c*d*(c^{-1})*d)^2=a*c*d*a*c*(d^{-1})*c*a*(d^{-1})=1  $ \\
\hline

\end{tabular}

\begin{tabular}{|l|c|} \hline
  \multicolumn{2}{|c|}{Table 3b} \\ \hline
  \multicolumn{2}{|c|}{Groups of order $32p$ without a normal sylow
$p$-subgroup}  \\ \hline
\multicolumn{2}{|c|}{II. Cases without a normal sylow subgroup (order 96
cases)}\\  \hline
\multicolumn{2}{|c|}{A. Cases arising as an extension from the group $A_4$
(7 cases)}\\
\multicolumn{2}{|c|}{[all cases with a $C_2$ action on the group $A_4$]}\\
\hline
Group & Automorphism Group \\ \hline
        20. $ A_4 @ 1^3 = S_4\times C_2\times C_2$&$    S_4\times S_4   $
\\
                                               &
\\
\hline

\end{tabular}\\

\begin{tabular}{|l|c|} \hline
  \multicolumn{2}{|c|}{Notes for Table 3b} \\ \hline
  \multicolumn{2}{|c|}{Relations for the non-direct-product cases, and some
comments}  \\
\multicolumn{2}{|c|}{on distinguishing these groups from one another.}  \\
\hline
\multicolumn{2}{|c|}{IIA. Here all actions on $A_4$ are by an operator of 
order 2. }\\  \hline
\multicolumn{2}{|c|}{$A_4: x^2=y^3=(y*x)^3=1;  C_2$ actions: $(a,x)=y^a*y=1
$}\\ \hline
\multicolumn{2}{|l|}{IIB. The relations here are: }\\

        \#27. [SL$(2,3) @ C_2] \times C_2 $ &$
a^3=a*b*a*((b*a*b)^{-1})=c^2=d^2= $     \\
                                   &$(c,d)=a^c*a*b*(a^{-1})=b^c*b=   $ \\
\hline
  \end{tabular}\\

\begin{tabular}{|l|c|} \hline
  \multicolumn{2}{|c|}{Table 3c} \\ \hline
  \multicolumn{2}{|c|}{Groups of order $32p$ without a normal sylow
$p$-subgroup}  \\ \hline
\multicolumn{2}{|c|}{III. Cases without a normal sylow 2-subgroup (order 96
cases)}\\  \hline
Group & Automorphism Group \\ \hline
      A. $p=5$      & \\                              &                   \\
               1.\quad $ 1^4 @ C_5 \times  C_2 $ & order 960 (complete group) \\
                                       &                                 \\
\hline
\end{tabular}%
\\

$\ddagger$ The numbers in this column refer to the number of this group as it appears in the normal subgroup list of the group Aut(G32) in the lists produced by CAYLEY. These numbers differ from the corresponding ones obtained from the GAP runs. For example, for group number 35 of order 32, the order 256 group appears as number 172 in CAYLEY and as number 185 in GAP's list. Likewise the order 128 group is 136 (CAYLEY) and 182 (GAP). \\ 

$\dagger $ For number 11, GAP says 20 characteristic cases and 60 total normal subgroups. Also, the number of the normal subgroup associated with the \textquotedblleft group [a]" according to GAP is number 58 in its normal subgroup listing for the group Aut(32 \#11). The original CAYLEY run file is missing, so we cannot check to see if there was a transcription error from the original CAYLEY run.  \\

\begin{tabular}{|l|} \hline
\qquad \qquad \qquad \qquad \qquad Notes for Table 4 \\ \hline
                                                                         \\
      number = number of the order 32 group as it appears in the          \\
           \qquad \qquad Hall-Senior Tables  [8].
      \\
      L(n)   = number of normal subgroups contained in the automorphism   \\
      \qquad \qquad     group of the order 32 group.
     \\
      \# char.  is the number of normal subgroups in the automorphism
\\
          \qquad \qquad   group of the group of order 32 that are
characteristic subgroups.                \\
                                                                          \\

             If the comment \textquotedblleft not a normal subgroup" appears, it means that \\
             the correct order and group are not given for the Aut($g$) \\
             factor; i.e., this factor is not a normal subgroup of             \\
             the automorphism group of the order 32 group in              \\
                                                                          \\
      \qquad \qquad \qquad        $G$ = $ C_p @ G[32]$.
          \\
      * not run.                                                          \\
                                                                          \\
      x wrong order cases:                                                \\
       \qquad \qquad 3a.  The calculation gives order 256 instead of 512.
      \\
       \qquad \qquad 8c.  The calculation gives order 32 instead of 256. \\
        \qquad \qquad 16ab. The calculation here gives order 64 instead of 128.  \\
        \qquad \qquad 19a.  The calculation here gives order 32 instead of 64. \\
                                                                          \\
         In all of these cases, however, it does appear that a group of the    \\
            correct type appears in the subgroup lattice of the           \\
            automorphism group of the appropriate group of order 32.      \\
            It would seem that the method used to calculate these         \\
            factors from the order 32 automorphism group is missing       \\
            some of these subgroups.                                      \\
                    \\
         In many cases the number of the normal subgroup in the subgroup  \\
            lattice of the automorphism group of the group of order 32    \\
            corresponding to the correct invariant factor given in        \\
            Table 2 above is given in the comments column.                \\
\hline
\end{tabular}\\

\begin{tabular}{|c|c|c|c|} \hline
  \multicolumn{4}{|c|}{Table 5} \\ \hline
  \multicolumn{4}{|c|}{Numbers of groups of order $32p^2$ arising from
groups of }  \\
\multicolumn{4}{|c|}{order 32 with an order 16 group action on the group
$C_p \times  C_p $}  \\ \hline
image &prime & number & original order 32 \\
     &    &  of cases & groups \\ \hline
   $ C_{16}$  &  $p \equiv  1$ mod(16) &     ?     &$   C_{32}$ and $C_{16} 
\times
C_2$   \\
            &  $p \equiv  1$ mod(8)  &     2     &                          
\\
\hline
   $ C_8 \times C_2$ &  $p \equiv  1$ mod(8)  &    14     & $ (3,1^2) [2], (3,2) [2]$ \\
            &                &           &  (4,1)[3],  \#20  [2] \\
            &                &           &  $(4,8|8,8)$[2],         \\
            &                &           &  \#22       [3].         \\
\hline
   $C_4 \times C_4 $ &  $p \equiv  1$ mod(4)  &     6     &  $(2^2,1)$,  
(3,2) [2]
   \\
            &                &           &  \#18       \#19 [2].    \\
\hline

   $C_4 \text{Y} Q_2$  &  $p \equiv  1$ mod(4)   &    10     & $ C_4 
\text{Y} Q_2 \times C_2, $
      \\
            &                &           & $ D_4 \times C_4,  Q_2 \times
C_4,$      \\
            &                &           &  \#16, \#36, \#37,       \\
            &                &           &  \#38, \#39, \#40,       \\
            &                &           &  \#41.                   \\
\hline
   $<2,2|2>$  &  $p \equiv  1$ mod(4)  &     4     &  $<2,2|2> \times\, C_2,$  
\#19,
   \\
            &                &           &  \#20,  \#21             \\
\hline
     $ D_8$    &  $p \equiv  7$ mod(8)  &     6     & $ D_8 \times C_2,  
(4,8|2,2),
$   \\
            &  $p \equiv  1$ mod(8)  &           &  \#29, $ D_{16}, QD_{16}, 
Q_8. $
  \\ \hline
     $ QD_8 $  &  $p \equiv  3$ mod(8)  &     4     & $ QD_8 \times C_2,  
(4,8|2,2),
$  \\
            &  $p \equiv  1$ mod(8)  &           &   \#28, \#30            
  \\
\hline
      $Q_4 $   &  $p \equiv  7$ mod(8)  &     3     & $ Q_4 \times C_2,$ \#28, \#29 \\
            &  $p \equiv  1$ mod(8)  &           &                          
\\
\hline
\multicolumn{4}{|c|}{ Note that for the last three cases one gets the same number}  \\
\multicolumn{4}{|c|}{for both sets of primes.}  \\ \hline
\multicolumn{4}{|c|}{If there is more than one group of order $32p^2$ coming
from a}  \\
\multicolumn{4}{|c|}{given order 32 group, then the number in [ ] in the last
column   }  \\
\multicolumn{4}{|c|}{after the group indicates the expected number of non-
       }  \\
\multicolumn{4}{|c|}{isomorphic groups of order $32p^2$ coming from this
group of      }  \\
\multicolumn{4}{|c|}{ order 32. }  \\ \hline
\end{tabular}

\begin{thebibliography}{99}
\bibitem{9} Becker, W., \textquotedblleft Automorphism groups for the groups of the form
($C_p \times  C_p$) @ $X$, where $|X|$ has order 32 and has a $D_4$ or $Q_2$
action on the $p$-group, paper in preparation. \\
\bibitem{12} Besche, Hans Ulrich, Eick, Bettina, and O'Brien, E. A., 
\textquotedblleft A millennium project: constructing small groups" is just one of
several articles dealing with the constructing and the listing of the results of 
their calculations.\\
\bibitem{14} Becker, W., \textquotedblleft Questions arising from games one can play
with groups and automorphism groups".\\
\end{thebibliography}
\end{document}